\documentclass[reqno]{amsart}
\usepackage{latexsym}
\usepackage{amsbsy}
\usepackage{amssymb}
\usepackage{amsfonts}
\usepackage{amsmath}
\usepackage[latin1]{inputenc}%needed for bibliography
\usepackage{graphicx}
\usepackage{a4wide}
 \newtheorem{thm}{Theorem}[section]
 \newtheorem{definition}[thm]{Definition}
 \newtheorem{lem}[thm]{Lemma}
 \newtheorem{prop}[thm]{Proposition}
 \newtheorem{cor}[thm]{Corollary}
 \newtheorem{rem}[thm]{Remark}
 \newtheorem{example}[thm]{Example}

 \newcommand{\bthm}{\begin{thm}}
 \newcommand{\ethm}{\end{thm}}
 \newcommand{\bd}{\begin{defin}}
 \newcommand{\ed}{\end{defin}}
 \newcommand{\blem}{\begin{lem}}
 \newcommand{\elem}{\end{lem}}
 \newcommand{\bcor}{\begin{cor}}
 \newcommand{\ecor}{\end{cor}}
 \newcommand{\bprop}{\begin{prop}}
 \newcommand{\eprop}{\end{prop}}
 \newcommand{\brem}{\begin{rem} \rm}
 \newcommand{\erem}{\end{rem}}
 \newcommand{\bex}{\begin{ex} \rm}
 \newcommand{\eex}{\end{ex}}

 \newcommand{\beq}{\begin{equation}}
 \newcommand{\eeq}{\end{equation} }

 \newcommand{\bea}{\begin{eqnarray}}
 \newcommand{\eea}{\end{eqnarray}}

 \newcommand{\beas}{\begin{eqnarray*}}
 \newcommand{\eeas}{\end{eqnarray*}}

 \newcommand{\beqs}{\begin{equation*}}
 \newcommand{\eeqs}{\end{equation*}}

 \newcommand{\bi}{\begin{itemize}}
 \newcommand{\ei}{\end{itemize}}

 \newcommand{\ben}{\begin{enumerate}}
 \newcommand{\een}{\end{enumerate}}

 \newcommand{\ba}{\begin{array}}
 \newcommand{\ea}{\end{array}}

\newcommand{\rd}{\mathbb{R}^d}

\begin{document}

\title[Suppleness of the sheaf]{Suppleness of the sheaf of algebras of generalized functions on manifolds}

\author{Stevan Pilipovi\'c}
          \address{Department of Mathematics, University of Novi Sad, Trg Dositeja Obradovi\'ca 4, Novi Sad, Serbia \\
              Tel.: +381-21-4852860,
              Fax: +381-21-6350458
               }
          \email{stevan.pilipovic@dmi.uns.ac.rs}
          % \thanks{The work of J. Aleksi\' c is supported by Ministry
           %   of Science and technological development, Republic of Serbia,
            %  project 144016}
\author{Milica \v Zigi\'c}
           \address{Department of Mathematics, University of Novi Sad, Trg Dositeja Obradovi\'ca 4, Novi Sad, Serbia \\
              Tel.: +381-21-4852790,
              Fax: +381-21-6350458}
           \email{milica.zigic@dmi.uns.ac.rs}

\begin{abstract}
       We show that the sheaves of algebras of generalized functions
$\Omega\to \mathcal{G}(\Omega)$ and $\Omega\to
\mathcal{G}^{\infty}(\Omega)$, $\Omega$ are open sets in a manifold
$X$, are supple, contrary to the non-suppleness of the sheaf of
distributions.
\end{abstract}

%\marginlabel{Mislim da je  recenica "In particular, this property is
%       satisfied by entropy solutions to heterogeneous scalar
%       conservation laws." u abstractu suvisna}

\maketitle

%%%%%%%%%%%%%%%%%%%%%%%%%%%%%%%%%%%%%%%%%%%%%%%%%%%%%%%%%%%%%%%%%%
 \section{Introduction and definitions} \label{sec:intro}
%%%%%%%%%%%%%%%%%%%%%%%%%%%%%%%%%%%%%%%%%%%%%%%%%%%%%%%%%%%%%%%%%%

The aim of this paper is to give a complete answer to the question
concerning suppleness of sheaves of certain generalized function
algebras. This question is discussed in \cite{msd} and here it is
completely solved. Note that Bros and Iagolnitzer \cite{bros}
conjectured that the analytic singular support (analytic wavefront
set) for distributions is decomposable. Bengel and Schapira
\cite{ben} have studied this decomposition by considering Cousin's
problem with bounds in a tuboid. In \cite{eida} authors have studied
microlocal decomposition for ultradistributions and
ultradifferentiable functions. They used the Laubin decomposition of
delta distribution \cite{laub} for the proof in this setting. We
consider in this paper the algebra $\mathcal{G}$ of generalized
functions containing the Schwartz distributions space $\mathcal{D}'$
as a subspace so that all the linear operations on $\mathcal{D}'$
are preserved within $\mathcal{G}$. We refer to \cite{bia},
\cite{col}, \cite{co1}, \cite{gfks}, \cite{gkos} and \cite{ober 001}
for the theory of generalized function algebras and applications to
non-linear and linear problems with non-smooth coefficients. Such
algebras are also called Colombeau algebras, since he was the first
one who introduced and analyzed such algebras. The geometric theory
of algebras of generalized functions \cite{gkos} is further
developed in papers \cite{kosv}, \cite{ks2}, \cite{ksv},
\cite{ksv2}, \cite{ksv3}, \cite{sv}. In these papers applications to
general relativity show the strong impact of the new approach
developed by the authors through the analysis of PDE on manifolds
with singular metrics and, in particular, in Lie group analysis of
differential equations (see \cite{gfks}, \cite{gkos}, \cite{gksv},
\cite{ko}). A version of this theory, which is the object of the
present article, is initiated in \cite{dd}, \cite{ks}. The sheaf
properties of generalized function algebras are investigated in
\cite{dp}, \cite{msd}.

In this paper we are interested in an important sheaf property, the
suppleness. It is known that the sheaves of Schwartz distributions
$\Omega\to \mathcal{D}'(\Omega)$ and of smooth functions $\Omega\to
C^\infty(\Omega)$, where $\Omega$ varies through all open sets of a
manifold $X$, are not supple. The extensions of these sheaves
$\Omega\to \mathcal{G}(\Omega)$ and $\Omega\to
\mathcal{G}^\infty(\Omega)$,
 $\Omega$ are open sets of a manifold $X$, which
are actually sheaves of algebras of generalized functions, are
supple. The proof of this assertion is the subject of this paper.

\subsection{Generalized functions on $\mathbb{R}^d$}\label{0.1}

We recall the main definitions. %We refer to the cited monographs for these definitions.
Let $\Omega$ be an open set in $\mathbb{R}^{d}$ and
$\mathcal{E}(\Omega)$ be the  space of nets of smooth functions.
Then the set of moderate nets  $\mathcal{E}_{M}(\Omega)$,
respectively of negligible nets  $\mathcal{N}(\Omega),$  consists of
nets $(f_{\varepsilon })_{\varepsilon\in(0,1)}\in
\mathcal{E}(\Omega)$ with the properties
$$
(\forall K\subset\subset \Omega)\;(\forall n\in\mathbb{N})\;(\exists
a\in\mathbb{R})\;(\sup\limits_{x\in K}|f^{(n)}_{\varepsilon
}(x)|=O(\varepsilon^{a})),
$$%
\[
\mbox{ respectively, } \;(\forall K\subset\subset \Omega)\;(\forall
n\in\mathbb{N})\;(\forall b\in\mathbb{R})\;(\sup\limits_{x\in
K}|f^{(n)}_{\varepsilon}(x)|=O(\varepsilon^{b}))
\]
($O$ is the Landau symbol "big O" and $K\subset\subset\Omega$ means
that $K$ is compact in $\Omega$ or that $\bar{K}$ is compact in
$\Omega$.)
 Both spaces are algebras and the latter is an ideal of  the former.

 The algebra of generalized functions
$\mathcal{G}(\Omega)$ is defined as the quotient
$\mathcal{G}(\Omega)=\mathcal{E}_M(\Omega)/\mathcal{N}(\Omega).$
This is also a differential algebra. If the nets
$(f_\varepsilon)_\varepsilon$ consist of constant functions on
$\Omega $ (i.e. supremums over the compact set $K$ reduce to the
absolute value), then one obtains
 the corresponding spaces  $\mathcal{E}_M$ and
$\mathcal{N}_0.$ They are algebras, $\mathcal{N}_0$ is an ideal in
$\mathcal{E}_M$ and, as a quotient, one obtains the algebra of
generalized complex numbers
$\bar{\mathbb{C}}=\mathcal{E}_M/\mathcal{N}_0$ (or
$\bar{\mathbb{R}}$). It is a ring, not a field.

The embedding of the Schwartz distributions in
$\mathcal{E}^{\prime}(\Omega)$ is realized through the sheaf
homomorphism $ \mathcal{E}^{\prime}(\Omega)\ni f\mapsto
[(f\ast\phi_{\varepsilon}|_{\Omega})_\varepsilon]\in\mathcal{G}(\Omega),
$ where the fixed  net of mollifiers
$(\phi_{\varepsilon})_{\varepsilon}$ is defined by
$\phi_{\varepsilon}=\varepsilon ^{-d}\phi(\cdot/\varepsilon),\;
\varepsilon<1,$ where $\phi\in\mathcal{S}(\mathbb{R}^{d})$ satisfies
$$\int
\phi(t)dt=1,\;\int t^{m}\phi(t)dt=0,m\in\mathbb{N}_{0}^{n},|m|>0.$$
($t^{m}=t_{1}^{m_{1}}...t_{n}^{m_{n}}$ and $|m|=m_{1}+...+m_{n}.$)
In fact $ \mathcal{E}^{\prime}(\Omega)$ is embedded into the space
$\mathcal{G}_c(\Omega)$ of compactly supported generalized
functions. This
 sheaf homomorphism,  extended onto $\mathcal{D}^{\prime}$, gives the embedding of $\mathcal{D}^{\prime
}(\Omega)$ into $\mathcal{G}(\Omega).$

The algebra of generalized functions $\mathcal{G}^\infty(\Omega)$ is
defined in \cite{ober 001} as the quotient of
$\mathcal{E}^\infty_{M}(\Omega)$ and $\mathcal{N}(\Omega),$ where
$\mathcal{E}^\infty_{M}(\Omega)$ consists of nets $(f_{\varepsilon
})_{\varepsilon\in(0,1)}\in \mathcal{E}(\Omega)^{(0,1)}$ with the
properties
$$
(\forall K\subset\subset \Omega) (\exists a\in\mathbb{R}) (\forall
n\in\mathbb{N}) (\sup_{x\in
K}|f^{(n)}_{\varepsilon}(x)|=O(\varepsilon^{a})),
$$

%this mapping is extended on ${\cal D}'(\Omega).$
Note that $\mathcal{G}^\infty$ is a subsheaf of $\mathcal{G}.$

\subsection{Generalized functions on a manifold}\label{0.2}

We will recall the main definitions and assertions following
\cite{gkos}. Let $X$ be a smooth Hausdorff paracompact manifold. We
denote by $\mathcal{U}=\{(V_{\alpha},\psi_{\alpha}):\alpha\in
\Lambda\}$ an atlas on $X$, $\Lambda$ is the index set.

We use $\mathcal{P}(X,E)$ to denote the space of linear differential
operators $\Gamma(X,E)\rightarrow\Gamma(X,E)$, where $E$ is a vector
bundle on $X$ and $\Gamma(X,E)$ is the space of smooth sections of
the vector bundle $E$ over $X$. Particularly, if $E=X\times
\mathbb{R}$ we write $\mathcal{P}(X)$ instead of $\mathcal{P}(X,E)$.
We denote by $\frak X(X)$ the space of smooth vector fields on $X$.

Let $\mathcal{E}(X):=(C^{\infty}(X))^{(0,1)}$ and
$(u_{\varepsilon})_{\varepsilon}\in \mathcal{E}(X)$. Then the
following statements are equivalent:
\begin{enumerate}
\item $(\forall K\subset\subset X)\;(\forall P\in
\mathcal{P}(X))\;(\exists N\in \mathbb{N})\;(\sup\limits_{p\in
K}|Pu_{\varepsilon}(p)|=O(\varepsilon^{-N}))$;

\item $(\forall K\subset\subset X)\;(\forall k\in \mathbb{N}_{0})\;(\exists N\in
\mathbb{N})\;(\forall\xi_{1},...,\xi_{k}\in
\frak{X}(X))$\\
$(\sup\limits_{p\in
K}|L_{\xi_{1}}...L_{\xi_{k}}u_{\varepsilon}(p)|=O(\varepsilon^{-N})),$
($L_{\xi_{i}}$ is the Lie derivative);

\item For any chart $(V,\psi)$: $(u_{\varepsilon}\circ
\psi^{-1})_{\varepsilon}\in \mathcal{E}_{M}(\psi(V))$.
\end{enumerate}

Denote by $\mathcal{E}_{M}(X)$ the subset of $\mathcal{E}(X)$
defined by any of the conditions 1, 2 or 3. We call it the space of
moderate nets on the manifold $X$. The space of negligible nets is
defined as:

$$\mathcal{N}(X):=\{(u_{\varepsilon})_{\varepsilon}\in \mathcal{E}_{M}(X):\forall K\subset\subset X\;\forall m\in \mathbb{N}:
\;\sup\limits_{p\in K}|u_{\varepsilon}(p)|=O(\varepsilon^{m})\}.$$

An algebra of generalized functions on the manifold $X$ is defined
as the quotient space
$\mathcal{G}(X):=\mathcal{E}_{M}(X)/\mathcal{N}(X).$ Elements of
$\mathcal{G}(X)$ are written as
$u=[(u_{\varepsilon})_{\varepsilon}]=(u_{\varepsilon})_{\varepsilon}+\mathcal{N}(X)$.
As one can expect, $\mathcal{E}_{M}(X)$ is a differential algebra
(with respect to Lie derivatives) and $\mathcal{N}(X)$ is a
differential ideal in it. Moreover, $\mathcal{E}_{M}(X)$ and
$\mathcal{N}(X)$ are invariant with respect to any $P\in
\mathcal{P}(X)$. Thus $Pu:=[(Pu_{\varepsilon})_{\varepsilon}]$ is a
well-defined element of $\mathcal{G}(X)$.

Let $u\in \mathcal{G}(X)$ and let $X'$ be an open set on a manifold
$X$. The restriction of a generalized function $u$, denoted by
$u|_{X'}\in \mathcal{G}(X')$, is represented by
$(u_{\varepsilon}|_{X'})_{\varepsilon}+\mathcal{N}(X')$. The support
of a generalized function $u$, denoted by ${\rm supp}\;u$, is
defined  as the complement of the union of open sets $X'\subseteq X$
such that $u|_{X'}=0$.

The algebra $\mathcal{G}^{\infty}(X)$ is defined as a subalgebra of
$\mathcal{G}(X)$ satisfying $u\in \mathcal{G}^{\infty}(X)$ if there
exists a representative $(u_\varepsilon)_\varepsilon$ of $u$ so that
for any chart $(U,\varphi)$, $(u_\varepsilon\circ
\varphi^{-1})_\varepsilon\in \mathcal{G}^{\infty}(\varphi(U))$.

Now we recall the sheaf properties of the space $\mathcal{G}(X)$
(see \cite{gkos}) and $\mathcal{G}^{\infty}(X)$.

A generalized function $u$ on $X$ allows the following local
description via the correspondence: $\mathcal{G}(X)\ni u\mapsto
(u_{\alpha})_{\alpha \in A}$, where $u_{\alpha}:=u\circ
\psi_{\alpha}^{-1}\in \mathcal{G}(\psi_{\alpha}(V_{\alpha}))$. We
call $u_{\alpha}$ the local expression of $u$ with respect to the
chart $(V_{\alpha},\psi_{\alpha})$. Then $\mathcal{G}(X)$ can be
identified with the set of all families $(u_{\alpha})_{\alpha}$ of
generalized functions $u_{\alpha}\in
\mathcal{G}(\psi_{\alpha}(V_{\alpha}))$ satisfying the
transformation law
$$u_{\alpha}|_{\psi_{\alpha}(V_{\alpha}\cap V_{\beta})}=
u_{\beta}\circ\psi_{\beta}\circ\psi_{\alpha}^{-1}|_{\psi_{\alpha}(V_{\alpha}\cap
V_{\beta})}$$ for all $\alpha,\beta\in A$ with $V_{\alpha}\cap
V_{\beta}\neq\emptyset$.

It is well known that $\Omega \to \mathcal{G}(\Omega)$, $\Omega$ are
open sets in $X$, is a fine and soft sheaf of $\mathbb{K}$-algebras
on $X$. Thus, $\mathcal{G}$ is defined directly as a quotient sheaf
of the sheaves of moderate modulo negligible sections. Similarly,
$\Omega\to \mathcal{G}^{\infty}(\Omega)$, $\Omega$ open in $X$, is a
fine and soft sheaf.

%%%%%%%%%%%%%%%%%%%%%%%%%%%%%%%%%%%%%%%%%%%%%%%%%%%%%%%%%%%%%%%%%%
 \section{Supple sheaves} \label{sec:supple}
%%%%%%%%%%%%%%%%%%%%%%%%%%%%%%%%%%%%%%%%%%%%%%%%%%%%%%%%%%%%%%%%%%

Recall \cite{war}, if $\mathcal{F}$ is a sheaf over the differential
manifold $X$ and $U\subset X$ open than a continuous map $f:U\to
\mathcal{F}$ such that $\pi\circ f=id$ is called a section of
$\mathcal{F}$ over $U$. The set of sections of $\mathcal{F}$ over
$U$ is denoted by $\Gamma(U,\mathcal{F})$.
\begin{definition}\label{savitljiv} Let $\mathcal{F}$ be a sheaf over the topological space $X$. Then, $\mathcal{F}$
is a supple sheaf if for all $f\in \Gamma(U,\mathcal{F})$, $U$ open
in $X$, the following is true: If ${\rm supp}\;f=Z=Z_{1}\cup Z_{2}$,
where $Z_{1}$ and $Z_{2}$ are arbitrary closed sets of $X$, then
there exist $f_{1},f_{2}\in \Gamma(U,\mathcal{F})$ such that $${\rm
supp}\;f_{1}\subseteq Z_{1},\; {\rm supp}\;f_{2}\subseteq
Z_{2}\;\mbox{ and}\; f=f_{1}+f_{2}.$$
\end{definition}

It will be shown that the sheaf of algebras $\Omega\to
\mathcal{G}(\Omega)$, $\Omega$ varies over all open sets of the
manifold $X$, is supple, but it is not flabby. It is well known that
$\mathcal{D}'$ is not supple. We give an example which shows this.

\begin{example}\label{dist koja nije supple}
Consider
$$f(x)=\sum\limits_{n=1}^{\infty}(\delta(x+\frac{1}{n^2})-\delta(x-\frac{1}{n^2}))+\delta(x),$$
where $\delta$ is the delta distribution. The support of this
distribution is
$$ Z=\{\frac{1}{n^2}:\;n\in \mathbb{N}\}\cup\{-\frac{1}{n^2}:\;n\in
\mathbb{N}\}\cup \{0\}.$$ One can see that the closed set  $Z$ is a
union of two closed sets $$\displaystyle
Z_{1}=\{\frac{1}{n^2}:\;n\in \mathbb{N}\}\cup\{0\}\;\mbox{ and}\;
\displaystyle Z_{2}=\{-\frac{1}{n^2}:\;n\in \mathbb{N}\}\cup\{0\}.$$
Then distributions $f_{1}$ and $f_{2}$ (in order to satisfy
Definition \ref{savitljiv}) should be of the form $$
f_{1}=\sum\limits_{n=1}^{\infty}\delta(x-\frac{1}{n^2})+C_{1}\delta(x)\;\mbox
{and} \;
f_{1}=\sum\limits_{n=1}^{\infty}\delta(x+\frac{1}{n^2})+C_{2}\delta(x),$$
since ${\rm supp\;}f_{1}\subseteq Z_{1}$, ${\rm
supp\;}f_{2}\subseteq Z_{2}$ and $f=f_{1}+f_{2}$. However, it is
known that $f_{1}$ and $f_{2}$ are not distributions since they are
infinite sums of shifted delta distributions so that their supports
have zero as the accumulation point.
\end{example}

We will prove the next theorem:

\begin{thm}\label{GX je savitljiv snop}$\Omega\to\mathcal{G}(\Omega)$, $\Omega$ open in $X$, is a supple sheaf.
\end{thm}

\textbf{Proof.} For the set $A$, we will denote by $A^\varepsilon$
the set $A^{\varepsilon}=\{x\in \rd
:\;d(x,A)<\varepsilon\},\;\varepsilon<1$, where $d$ is a distance on
$X$. The notation $L(x,r)$ stands for the open ball of radius $r>0$
centered in $x\in X$, i.e. $L(x,r)=\{y\in X:\;d(x,y)<r\}$.

We divide the proof of this theorem into two parts (I) and (II) and
use the following two simple assertions:
\begin{enumerate}
    \item Let $A$ be a measurable set in $\rd$. Then there exists a
    generalized function $\eta=[(\eta_\varepsilon)_\varepsilon]\in
    \mathcal{G}(\rd)$ such that
    $$\eta_{\varepsilon}(x):=\left\{\begin{array}{ll}
    1, & x\in A\\
    0, & x\in \rd\setminus A^{\varepsilon}
    \end{array}\right.$$
and $0\leq|\eta_\varepsilon(x)|\leq 1,\;x\in \rd$. More precisely,
$\eta_\varepsilon$ is defined to be $1_A\ast \phi_\varepsilon$,
where $1_A$ is the characteristic function of $A$, $\phi$ is a
compactly supported smooth function so that $\int_{\rd}\phi(t)dt=1$
and $\phi_\varepsilon(x)=1/\varepsilon^d\phi(x/\varepsilon)$.

\item
    Let $\delta>0$ and $Z_{1}$ and $Z_{2}$ arbitrary closed sets of
    $\rd$. Then there exists a closed set
    $\widetilde{Z_{1}^{\delta}}\supset Z_{1}$ such that $Z_{1}\cap
    Z_{2}=\widetilde{Z_{1}^{\delta}}\cap Z_{2}$ and $d(x,Z_{1})\leq
    \delta$, $x\in \widetilde{Z_{1}^{\delta}}$.

 We define $$\widetilde{Z_{1}^{\delta}}=\{x\in \rd
:\;d(x,Z_{1})\leq \delta\;\wedge d(x,Z_{1})\leq d(x,Z_{2})\}.$$

\end{enumerate}

     (I)   Now we show that $\Omega\to\mathcal{G}(\Omega)$, $\Omega$ open
        in $\rd$, is a supple sheaf.

        It is enough to prove the assertion for $U=\rd$ and $f\in \mathcal{G}(\rd)=\Gamma(\rd,\mathcal{G})$ with the property ${\rm
        supp}\;f=Z$ and let $Z=Z_{1}\cup Z_{2}$, where $Z_{1}$ and $Z_{2}$
        are arbitrary closed sets. Let $\delta>0$ and define
        $\widetilde{Z_{1}^{\delta}}$ as in Assertion 2. Next by Assertion 1,
        let $\eta_{\varepsilon}\in C^{\infty}(\rd),\;\varepsilon\in(0,1)$,
        such that
        $$\eta_{\varepsilon}(x)=\left\{\begin{array}{ll}
        1, & x\in \widetilde{Z_{1}^{\delta}}\\
        0, & x\in X\setminus
        (\widetilde{Z_{1}^{\delta}})^{\varepsilon},\;\varepsilon\in (0,1),
        \end{array}\right.$$
        with $(\eta_\varepsilon)_\varepsilon\in \mathcal{E}_M(\rd)$.

        Let $f_{1}=[(f_{\varepsilon}\eta_{\varepsilon})_{\varepsilon}]$ and
        $f_{2}=[(f_{\varepsilon}(1-\eta_{\varepsilon}))_{\varepsilon}]$.
        Then $f_{1},f_{2}\in \mathcal{G}(\rd)$ and
        $f=[(f_{\varepsilon}\eta_{\varepsilon})_{\varepsilon}]+
        [(f_{\varepsilon}(1-\eta_{\varepsilon}))_{\varepsilon}]$. So, we
        have to show that ${\rm supp}\;(f_{1})\subseteq Z_{1}$ and ${\rm
        supp}\;(f_{2})\subseteq Z_{2}$.

        We show the inclusion ${\rm supp}\;(f_{1})\subseteq Z_{1}$ by
        showing that for any point $x\notin Z_{1}$ there exists a
        neighborhood $X'$ of $x$ such that
        $(f_{1\varepsilon}|_{X'})_{\varepsilon}\in \mathcal{N}(X')$
        (according to the definition of the support this means $x\notin {\rm
        supp}\;f_{1}$). Let $x\in \rd \setminus Z_{1}$. Then we have
        $A=d(x,Z_{1})>0$ and there are two possibilities: $x\in
        \widetilde{Z_{1}^{\delta}}$ and $x\notin
        \widetilde{Z_{1}^{\delta}}$. If $x\in \widetilde{Z_{1}^{\delta}}$
        then the ball $X'=L(x,\frac{A}{2})$ has no intersection with $Z_{1}$
        and $Z_{2}$ because in the set $\widetilde{Z_{1}^{\delta}}$ we have
         $d(x,Z_{1})\leq d(x,Z_{2})$. So $X'$ does not intersect the set
        $Z=Z_{1}\cup Z_{2}$. From
        $$|f_{1\varepsilon}(y)|=|f_{\varepsilon}(y)\eta_{\varepsilon}(y)|\leq
        |f_{\varepsilon}(y)| \;\mbox{for all}\; y\in X'$$ and
        $(f_{\varepsilon}|_{X'})_{\varepsilon}\in \mathcal{N}(X')$, we have
        $(f_{1\varepsilon}|_{X'})_{\varepsilon}\in \mathcal{N}(X')$.
        Finally, if $x\notin\widetilde{Z_{1}^{\delta}}$ then
        $B=d(x,\widetilde{Z_{1}^{\delta}})>0$ (since
        $\widetilde{Z_{1}^{\delta}}$ is a closed set). Let
        $X'=L(x,\frac{B}{2})$. Then
        $$f_{1\varepsilon}(y)=f_{\varepsilon}(y)\eta_{\varepsilon}(y)=f_{\varepsilon}(y)\cdot
        0=0,\;y\in X',$$ where $\varepsilon<\frac{B}{2}$. Again, we have
        $(f_{1\varepsilon}|_{X'})_{\varepsilon}\in \mathcal{N}(X')$.

        Similarly, we show the second inclusion ${\rm
        supp}\;(f_{2})\subseteq Z_{2}$. Let $x\notin Z_{2}$. Our aim is to
        show that $x\notin {\rm supp}\;f_{2}$. There are two possibilities:
        $x\notin Z_{1}$ and $x\in Z_{1}$. If
        $x\notin Z_{1}$, we also have $x\notin Z_{2}$ and so $x\notin Z$.
        Then there exists a neighborhood $W$ of $x$ such that
        $(f_{\varepsilon}|_{W})_{\varepsilon}\in \mathcal{N}(W)$, since
        ${\rm supp}\;f\subseteq Z$. We take $X'=W$. Then, clearly,
        $$|f_{2\varepsilon}(y)|=|f_{\varepsilon}(y)(1-\eta_{\varepsilon}(y))|\leq
        |f_{\varepsilon}(y)|,\;y\in X',\;\varepsilon<1.$$ Since
        $(f_{\varepsilon}|_{X'})_{\varepsilon}\in \mathcal{N}(X')$ we also
        have that $(f_{2\varepsilon}|_{X'})_{\varepsilon}\in
        \mathcal{N}(X')$. The second possibility is $x\in Z_{1}$. Let
        $H=d(x,Z_{2})>0$ and $H'=\min \{H,\delta\}$. Note that
        $d(x,Z_{1})=0$ since $x\in Z_{1}$. Then for
        $X'=L(x,\frac{H'}{2})$ we have $X'\subseteq
        \widetilde{Z_{1}^{\delta}}$ (since $H'\leq \delta$ and for all $y\in
        X'$ holds $d(y,Z_{1})<\frac{H'}{2}\leq d(y,Z_{2})$). So $X'$ has no
        intersection with $Z_{2}$ (since $H'\leq H$). We have
        $$f_{2\varepsilon}(y)=f_{\varepsilon}(y)(1-\eta_{\varepsilon}(y))=f_{\varepsilon}(y)\cdot
        0=0,\;y\in X',\;\varepsilon<1.$$ Again,
        $(f_{2\varepsilon}|_{X'})_{\varepsilon}\in \mathcal{N}(X')$. This
        finishes the proof of the suppleness of $\Omega
        \to\mathcal{G}(\Omega)$, $\Omega$ open in $\rd$.

(II) Let $Z\subseteq X$ be closed. Let $Z=Z_1\cup Z_2$ where $Z_1$
and $Z_2$ are closed and let $f\in
\mathcal{G}(X)=\Gamma(X,\mathcal{G})$ such that ${\rm
supp}\;(f)\subseteq Z$. Cover $X$ by a family of chart neighborhoods
$\mathcal{U}$ and let $\{\tilde{\chi}_{\alpha}:\;\alpha\in
\Lambda\}$ be a partition of unity subordinated to $\mathcal{U}$.
Set
$$\chi_\alpha=\frac{\tilde{\chi}_\alpha}{(\sum_{\alpha\in \Lambda}\tilde{\chi}^{2}_{\alpha})^{1/2}}.$$
So, we obtain the family of functions $\{\chi_\alpha:\;\alpha\in
\Lambda\}$ such that $\{{\rm supp}\;(\chi_\alpha):\;\alpha\in
\Lambda\}$ is locally finite and $\sum_{\alpha\in
\Lambda}\chi^{2}_{\alpha}=1$. Hence, one can write
\begin{equation}\label{rastavljane f} f=\sum_{\alpha\in \Lambda}
\chi^{2}_{\alpha} f=\sum_{\alpha\in \Lambda}\chi_\alpha(\chi_\alpha
f).\end{equation} For the functions $\chi_\alpha f\;\alpha\in
\Lambda$ we have ${\rm supp}(\chi_\alpha f)$ is closed in some
$U_\alpha\in \mathcal{U}$, where the results hold (due to part (I)
of this proof). Precisely, one can see ${\rm supp}(\chi_\alpha
f)\subseteq Z\cap U_\alpha$ as
$${\rm supp}(\chi_\alpha f)=(Z_1\cap ({\rm supp}(\chi_\alpha
f)))\cup (Z_2\cap ({\rm supp}(\chi_\alpha f)))\subseteq(Z_1\cap
U_\alpha)\cup(Z_2\cap U_\alpha),$$ where the sets $Z_1\cap ({\rm
supp}(\chi_\alpha f))$ and $Z_2\cap ({\rm supp}(\chi_\alpha f))$ are
closed.

Applying part (I) of this proof to $\chi_\alpha f,\;\alpha\in
\Lambda$ we obtain $f^{\alpha}_{1},f^{\alpha}_{2}\in
\mathcal{G}(U_\alpha)$ such that
$$\chi_\alpha f=f^{\alpha}_{1}+f^{\alpha}_{2},\;{\rm
supp}(f^{\alpha}_{1})\subseteq Z_1\cap U_\alpha\subseteq Z_1,\;{\rm
supp}(f^{\alpha}_{2})\subseteq Z_2\cap U_\alpha\subseteq Z_2.$$
According to (\ref{rastavljane f})

$f=\sum_{\alpha\in
\Lambda}\chi_\alpha(f^{\alpha}_{1}+f^{\alpha}_{2})=\sum_{\alpha\in
\Lambda}\chi_\alpha f^{\alpha}_{1}+\sum_{\alpha\in \Lambda}
\chi_\alpha f^{\alpha}_{2}$.

Set $f_1=\sum_{\alpha\in \Lambda}\chi_\alpha f^{\alpha}_{1}$ and
$f_2=\sum_{\alpha\in \Lambda} \chi_\alpha f^{\alpha}_{2}$. Then
$f_1,f_2\in \mathcal{G}(X)$ and ${\rm supp}f_1\subseteq Z_1,\;{\rm
supp}f_2\subseteq Z_2$.

$\hfill\Box$

\begin{thm}\label{GXinf je savitljiv snop}$\Omega\to\mathcal{G}^{\infty}(\Omega)$, $\Omega$ open in $X$, is a supple sheaf.
\end{thm}

\textbf{Proof.} Suppleness of the sheaf $\Omega\to
\mathcal{G}^{\infty}(\Omega)$ can be proved using the same ideas as
in the proof of Theorem \ref{GX je savitljiv snop}. Let $Z\subseteq
X$ be closed. Let $Z=Z_1\cup Z_2$ where $Z_1$ and $Z_2$ are closed
and let $f\in \mathcal{G}^\infty(X)$ such that ${\rm
supp}\;(f)\subseteq Z$. Now we have to construct generalized
functions $f_1, f_2\in \mathcal{G}^\infty(X)$ such that ${\rm
supp}f_1\subseteq Z_1,\;{\rm supp}f_2\subseteq Z_2$. In order to
obtain $f_1, f_2\in \mathcal{G}^\infty(X)$ we will take $\hat{\eta}$
to be a generalized function from $ \mathcal{G}^{\infty}(X)$ (see
Assertion 1). We will
 replace $\varepsilon$ by $|\ln\varepsilon|^{-1}$ and then the
generalized function
$\hat{\eta}=[(\hat{\eta}_\varepsilon)_\varepsilon]$ will be
$$\hat{\eta}_{\varepsilon}(x)=\left\{\begin{array}{ll}
1, & x\in A\\
0, & x\in \rd\setminus A_{|\ln\varepsilon|^{-1}}
\end{array}\right..$$

This finishes the proof, since for
$f=[(f_\varepsilon)_\varepsilon]\in\mathcal{G}^\infty(X)$
generalized functions
$\hat{f}_{1}=[(f_{\varepsilon}\hat{\eta}_{\varepsilon})_{\varepsilon}]$
and
$\hat{f}_{2}=[(f_{\varepsilon}(1-\hat{\eta}_{\varepsilon}))_{\varepsilon}]$
are in $\mathcal{G}^\infty(X)$ and following the proof of Theorem
\ref{GX je savitljiv snop} (replacing $\varepsilon$ by
$|\ln\varepsilon|^{-1}$ ) we obtain ${\rm
supp}\;\hat{f}_{1}\subseteq Z_{1}$ and ${\rm
supp}\;\hat{f}_{2}\subseteq Z_{2}$.

$\hfill\Box$

Let us remark at the end that $\mathcal{G}(X)$ and
$\mathcal{G}^{\infty}(X)$ are not flabby sheaves (see \cite{msd},
Remark on page 95). If we take $X=\mathbb{R}$ and $X'=(0,\infty)$
then one can not extend the generalized function
$[(\varepsilon^{-1/x})_{\varepsilon}]$, defined on $(0,\infty)$, to
the whole space $\mathbb{R}$.

%%%%%%%%%%%%%%%%%%%%%%%%%%%%%%%%%%%%%%%%%%%%%%%%%%%%%%%%%%%%%%%%%%
 \subsection*{Acknowledgement}
%%%%%%%%%%%%%%%%%%%%%%%%%%%%%%%%%%%%%%%%%%%%%%%%%%%%%%%%%%%%%%%%%%

 The research is supported by Ministry
of Science and Technological Development, Republic of Serbia,
project 144016.

%%%%%%%%%%%%%%%%%%%%%%%%%%%%%%%%%%%%%%%%%%%%%%%%%%%%%%%%%%%%%%%%%%
%%%%%%%%%%%%%%%%%%%%%%%%%%%%%%%%%%%%%%%%%%%%%%%%%%%%%%%%%%%%%%%%%%

%%%%%%%%%%%%%%%%%%%%%%%%%%%%%%%%%%%%%%%%%%%%%%%%%%%%%%%%%%%%%%%%%%
%%%%%%%%%%%%%%%%%%%%%%%%%%%%%%%%%%%%%%%%%%%%%%%%%%%%%%%%%%%%%%%%%%

\begin{thebibliography}{00}

\bibitem{ben}
G. Bengel, P. Schapira, D\'ecomposition microlocale analytique des
distributions, Ann. Inst. Fourier, Grenoble 29, 101-124 (1979)

\bibitem{bia}
H. A. Biagioni, A Nonlinear Theory of Generalized Functions,
Springer-Verlag, Berlin-Hedelberg-New York (1990)

\bibitem{bros}
J. Bros, D. Iagolnitzer, Support essentiel et structure analytique
des distribution, S\'eminaire Goulaouic-Lions-Schwartz, Expos\'e 18
(1975-1976)

\bibitem{col}
J. F. Colombeau, New Generalized Functions and Multiplications of
Distributions, North Holland, Amsterdam (1984)

\bibitem{co1}
J. F. Colombeau, Elementary Introduction in New Generalized
Functions, North Holland, Amsterdam (1985)

\bibitem{dd}
J.W. De Roever, M. Damsma,  Colombeau algebras on a
$C^\infty$-manifold, Indag. Math. N.S. 2, 341-358 (1991)

\bibitem{dp}
N. Djapi\'c, S. Pilipovi\'c,  Microlocal analysis of Colombeau's
generalized functions on a manifold, Indag. Math. N.S. 7, 293--309
(1996)

\bibitem{eida}
A. Eida, S. Pilipovi\'c, On the microlocal decomposition of same
classes of hyperfunctions, Math. Proc. Camb. Phil. Soc. 125, 455-461
(1999)

\bibitem{gfks}
M. Grosser, E. Farkas, M. Kunzinger, R. Steinbauer,  On the
foundations of nonlinear generalized functions I, II, Mem. Amer.
Math. Soc. 153 (2001)

\bibitem{gkos}
M. Grosser, M. Kunzinger, M. Oberguggenberger, R. Steinbauer,
Geometric Theory of Generalized Functions with Applications to
General Relativity, Kluwer Academic Publishers, Dordrecht (2001)

\bibitem{gksv}
M. Grosser, M. Kunzinger, R. Steinbauer, J. Vickers,  A global
theory of algebras of generalized functions, Adv. Math. 166, 50-72
(2002)

\bibitem{ko}
M. Kunzinger, M. Oberguggenberger,  Group analysis of differential
equations and generalized functions, SIAM J. Math. Anal. 31,
1192-1213 (2000)

\bibitem{kosv}
M. Kunzinger, M. Oberguggenberger, R. Steinbauer, J. Vickers,
Generalized flows and singular ODEs on differentiable manifolds,
Acta Appl. Math. 80, 221-241 (2004)

\bibitem{ks}
M. Kunzinger, R. Steinbauer, Foundations of a nonlinear
distributional geometry, Acta Appl. Math. 71, 179-206 (2002)

\bibitem{ks2}
M. Kunzinger, R. Steinbauer,  Generalized pseudo-Riemannian
geometry, Trans. Amer. Math. Soc. 354, 4179-4199 (2002)

\bibitem{ksv}
M. Kunzinger, R. Steinbauer, J. Vickers, Sheaves of nonlinear
generalized functions and manifold-valued distributions, Trans.
Amer. Math. Soc. 361, 5177-5192 (2009)

\bibitem{ksv2}
M. Kunzinger, R. Steinbauer, J. Vickers, Generalised connections and
curvature, Math. Proc. Cambridge Philos. Soc. 139, 497--521 (2005)

\bibitem{ksv3}
M. Kunzinger, R. Steinbauer, J. Vickers, Intrinsic characterization
of manifold-valued generalized functions, Proc. London Math. Soc.
87, 451-470 (2003)

\bibitem{laub}
P. Laubin, Front d'onde analytique et d\'ecomposition microlocale
des distributions, Ann. Inst. Fourier, Grenoble 33, 179-199 (1983)

\bibitem{ober 001}
M. Oberguggenberger, Multiplication of distributions and application
to partial differential equations, Pitman Res. Notes Math. Ser. 259,
Longman, Harlow (1992)

\bibitem{msd}
M. Oberguggenberger, S. Pilipovi\'c, D. Scarpalezos,  Local
properties of Colombeau generalized functions, Math. Nachr. 256,
88-99 (2003)

\bibitem{sv}
R. Steinbauer, J. Vickers, The use of generalized functions and
distributions in general relativity, Classical Quantum Gravity 23,
91-114 (2006)

\bibitem{war}
F. W. Warner, Foundations of Differentiable Manifolds and Lie
Groups, Scott, Foresman and Company, Glenview, Illinois, London
(1971)

\end{thebibliography}
 \end{document}